\documentclass{article}
\usepackage{amsfonts,amsmath,amsthm,amscd}
\usepackage{graphicx}

\def\R{{\mathbb R}}
\def\C{{\mathbb C}}

\def\Z{{\mathbb Z}}

\newcommand{\eps}{\varepsilon}
\newcommand{\pih}{\frac{\pi}{2}}

\begin{document}

\title{On the moment-angle manifolds with positive Ricci curvature}
\author{Ya.~V.~Bazaikin\thanks{Author was supported by the Leading
Scientific Schools grant NSh-7256.2010.1, SB RAS Integration grant
N46 and the RFBR (grants 09-01-00598-a and 10-01-92102-Yaf-a)},
I.~V.~Matvienko\thanks{Author was supported by the Leading
Scientific Schools grant NSh-7256.2010.1, SB RAS Integration grant
N46 and the RFBR (grant 09-01-00142-a)}}
\date{}
\maketitle

%\tableofcontents
%\newpage

\section{Introduction}

An interesting problem in Riemannian geometry is the question
of the topological complexity of the manifolds admitting Riemannian metrics
of positive Ricci curvature. For example, Sha and Yang
\cite{Sha-Yang} constructed metrics of positive Ricci curvature on
the connected sums of arbitrary number of copies of $S^n \times S^m$,
assuming $n, m \geq 2$ fixed. This implies the absence of a priori estimates for the Betti numbers. Wraith \cite{Wraith} generalized Sha and Yang’s examples and demonstrated the existence of Riemannian metrics with positive Ricci curvature on the connected sums
$\#_{i=1}^N S^{n_i} \times S^{m_i}$, where $n_i, m_i \geq 3$.

In the present article, we construct Riemannian metrics of positive Ricci
curvature on certain moment-angle manifolds. A moment-angle manifold
$Z_P$ is constructed on the basis of a polyhedron $P$, and there is the
canonical free action of the torus $T$ on $Z_P$ such that $Z_P/T=P$. The construction of the moment-angle manifolds is described in more detail in the next section. The main goal of the article is to establish the following result.

\vskip0.2cm

{\bf Theorem}. {\it Let $P$ be an octahedron obtained from a
three-dimensional cube by cutting off small neighborhoods of two
edges lying on skew lines. Then, the $11$-dimensional moment-angle
manifold $Z_P$ admits a Riemannian metric with positive Ricci
cur\-va\-tu\-re. }

\vskip0.2cm

Besides this example, we consider other moment-angle manifolds with
the metrics of positive Ricci curvature. However, the space
indicated in the theorem is of particular interest, because of its
topological complexity. For example, as is shown in \cite{Baskakov,
Buchshtaber-Panov}, it contains nontrivial Massey triple products in
cohomology. Hence, the constructed $Z_P$ is non-formal moment-angle
manifold.

In the conclusion of the introduction we would like to pose two
questions, which seem natural to us.

\vskip0.2cm

{\bf Question 1}. {\it Does there exist a metric of positive Ricci curvature on every moment-angle manifold?}

\vskip0.2cm

We believe that the answer to this question is positive. As one of the facts suggesting this (besides the results of the present article), we indicate the following. If $P$ is a two-dimensional polygon or $P$ is obtained from a multidimensional tetrahedron by multiple applying the operation of cutting off the neighborhood of a vertex, then, as is shown in \cite{Bosio-Meersseman},   $Z_P$ is diffeomorphic to a certain connected sum of the products of spheres of various dimensions, and the positive answer to the question follows from Wraith's
paper \cite{Wraith}.

    The concept of a quasitoric manifold \cite{Buchshtaber-Panov} is
closely related to the moment-angle manifolds. The second question is
more difficult and not so clear:

\vskip0.2cm

{\bf Question 2}. {\it Does there exist a metric of positive Ricci curvature
on every quasitoric manifold? }

\vskip0.2cm

To indicate that the positive answer is probable, we only point out that
the Riemannian metrics with positive Ricci curvature were constructed in
\cite{Sha-Yang2} on every simply connected four-dimensional
quasitoric manifold. Moreover, it was proved in
\cite{Bazaikin-Matvienko} that those metrics can be chosen to be
invariant with respect to an arbitrary prescribed action of $T^2$.

The authors are grateful to T.~E.~Panov for numerous useful discussions.

\section{Moment-angle manifolds}
The concept of a moment-angle manifold was introduced in
\cite{Davis-Januszkievicz}. The detailed treatment of the moment-angle
manifolds and the closely related quasitoric manifolds is given in
\cite{Buchshtaber-Panov}. Here we give only the necessary definitions and
properties.

Let $P$ be an $n$-dimensional polyhedron in the Euclidean space $\R^n$,
defined by the system of inequalities
$$
\sum_{i=1}^n a_{ij}x^j + b_i\geq 0, i=1, \ldots, m
$$
which are in general position. Let $P_1, \ldots, P_m$ be its faces of codimension one. For $p \in P$ let $G(p)$ denote the smallest (by
inclusion) face containing the point $p$. Consider $X=P \times T^m$,
where $T^m=\{ (z_1, \ldots, z_m) : z_i \in \C, |z_i|=1 \}$ is the
standard $m$-dimensional torus with the numeration of coordinates corresponding to the numeration of the codimension-one faces of $P$. Let
$T^{F_i}$ denote the circle $S^1 \subseteq T^m$ corresponding to the $i$-th
codimension-one face. Now, for each face $G$ of the polyhedron $P$ we put
$$
T^G=\prod_{G \subseteq F_i} T^{F_i} \subset T^m.
$$
We identify the points of $X$ in the following way:
$$
(p,z_i)\sim (p',z_i') \mbox{\ \, iff \ } p=p' \mbox{\ and \ }
\bar{z_i}z_i' \in T^{G(p)}.
$$
One can prove \cite{Buchshtaber-Panov, Davis-Januszkievicz} that
there is a canonical structure of a topological manifold on the
factor-space $Z_P=(P \times T^m)/\sim$, and the natural $T^m$ action
by shifts on itself induces a continuous action on $Z_P$. The
subgroup $T^G(p)$ is the stabilizer of the point $(p,Z) \in Z_P$
and, clearly,  $Z_P/T^m=P$. Moreover, $Z_P$ can be equipped with the
structure of a smooth manifold, such that the natural $T^m$ action
becomes smooth \cite{Buchshtaber-Panov}. The manifold $Z_P$ (with
some $T^m$-invariant smooth structure) is called a {\it moment-angle
manifold}.

Now, if we consider a torus $T^{m-n} \subset T^m$, which acts freely
on $Z_P$, then $M=Z_P /T^{m-n}$ will be a quasitoric manifold, and,
conversely, every quasitoric manifold can be constructed from a
certain moment-angle manifold \cite{ Davis-Januszkievicz}. In this
case, there arises a principal bundle with the structural group
$T^{m-n}$
$$
\pi: Z_P \rightarrow M.
$$
Analogously, if the factor-space $Z_P/T^{m-n}$ is an orbifold, then
we will call it a quasitoric orbifold. In this case, $\pi$ is a
principal toric bundle in the sense of orbifolds.

{\bf Example 0}. Let $P=I=[0,1]$ be a segment on the real line. In this case, we
can easily understand that $Z_P=S^3=\{ (u, v) \in \C^2 : |u|^2+|v|^2=1
\}$. The torus $T^2=\{ (z_1, z_2) \in \C^2 : |z_1|=|z_2|=1\}$ acts
on $S^3$ in the standard way
$$
(z_1,z_2) \in T^2: (u,v) \mapsto (u z_1, v z_2).
$$

{\bf Example 1}. Let $P=I^3$ be a three-dimensional cube. Then
$Z_P=S^3 \times S^3 \times S^3$. The torus $T^6=T^2 \times T^2 \times
T^2$ acts on $Z_P$. If we consider a circle $S^1$, diagonally embedded in $T^2$, then $S^1 \times S^1 \times S^1$ acts freely on $Z_P$, and we obtain a quasitoric manifold $M_1=S^2 \times S^2 \times S^2$.

The manifold $M_1$ can be modified in the following way. The torus
$T^3=T^6/(S^1 \times S^1 \times S^1)$ acts on $M_1$ (every circle of
the torus acts on its own sphere $S^2$ by rotations around the polar
axis). It is obvious that $M/T^3=P$. Consider a subgroup
$\Gamma=\Z_3$ of $T^3$ generated by the element $(\omega, \omega,
\omega) \in T^3$, where $\omega=e^{2 \pi i/3}$. Then,
$M'_1=M_1/\Gamma$ is a quasitoric orbifold with eight singular
points, whose neighborhoods look like $\C^3/\Z_3$. If we resolve
every singular point with the help of the "blow-up" operation (it
means that we cut out some neighborhood of a singular point and,
instead of it, paste in a cube of the canonical complex line bundle
over $\C P^2$  --- the details of this construction can be found in
the next section), then we arrive at a quasitoric manifold $N_1$.
Now, the polyhedron $Q_1=N_1/\Z_3$, associated with the quasitoric
manifold $N_1$, can be obtained from the cube $P$ by cutting off all
the vertices. We obtain a principal $T^3$-bundle
$$
\pi_1: Z_{Q_1} \rightarrow N_1.
$$

{\bf Example 2}. Let us consider a polyhedron $Q_2$, obtained from the
cube $P$ by cutting off the neighborhood of a certain vertex. To describe
$Z_{Q_2}$, we consider the following construction. Let the torus $T^3$
act on $Z_P=S^3 \times S^3 \times S^3$ as follows
$$
(z_1, z_2, z_3) \in T^3 : \left( \begin{array}{rr}
  u_1 & v_1 \\
  u_2 & v_2 \\
  u_3 & v_3 \\
\end{array} \right) \mapsto \left( \begin{array}{rr}
  z_1 z_2^2    \; u_1 &  \bar z_1 \; v_1 \\
  z_2 z_3^2    \; u_2 &  \bar z_2 \; v_2 \\
  \bar z_1 z_3 \; u_3 &  \bar z_3 \; v_3 \\
\end{array} \right),
$$
where the rows of matrices represent the coordinates of the corresponding spheres
$S^3$.

\vskip0.2cm

{\bf Lemma 1}. {\it The action described above is free everywhere in
$Z_P$ outside of the submanifold
$$
F = \left\{ \left( \begin{array}{rr}
  u_1 & 0 \\
  u_2 & 0 \\
  u_3 & 0 \\
\end{array} \right) : |u_1|=|u_2|=|u_3|=1 \right\}.
$$
The points of the submanifold $F$ have the stabilizer $Fix(F)= \Gamma =
\Z_3$ generated by the element $(\omega, \omega, \omega)$, where
$\omega = e^{2\pi i/3}$.}

\vskip0.2cm

Thus, the factor-space $M'_2=Z_P/T^3$ is a quasitoric orbifold with
one singular point $p=F/T^3$, whose neighborhood is diffeomorphic to
$\C^3/\Z_3$. If we resolve this singularity by blowing-up the
neighborhood of $p$, then we get a quasitoric manifold $N_2$
corresponding to $N_2/T^3=Q_2$ and obtain a principal $T^4$-bundle
$$
\pi_2: Z_{Q_2} \rightarrow N_2.
$$
In \cite{Bosio-Meersseman} (using the results of
\cite{Buchshtaber-Panov}) the topology of the space $Z_{Q_2}$ was
investigated. In particular, $Z_{Q_2}$ can not be represented as
connected sum of products of spheres. Further properties of the
topology of $Z_{Q_2}$ was studied in \cite{Gitler-Medrano}.

{\bf Example 3}. Let us consider a polyhedron $Q_3$, obtained from
the cube $P$ by cutting off the neighborhoods of two edges lying on
skew lines. To describe $Z_{Q_3}$, we consider the following
construction. Let the torus $T^3$ act on $Z_P=S^3 \times S^3 \times
S^3$ as follows
$$
(z_1, z_2, z_3) \in T^3_2 : \left( \begin{array}{rr}
  u_1 & v_1 \\
  u_2 & v_2 \\
  u_3 & v_3 \\
\end{array} \right) \mapsto \left( \begin{array}{rr}
  z_1               \; u_1 &  z_1 z_2 \bar{z}_3 \; v_1 \\
  z_2 \bar{z}_3     \; u_2 &  \bar{z}_1 z_2     \; v_2 \\
  \bar{z}_1 z_2 z_3 \; u_3 &  z_3               \; v_3 \\
\end{array} \right),\eqno{(1})
$$
where, as before, the rows of matrices represent the coordinates of the corresponding spheres
$S^3$.

\vskip0.2cm

{\bf Lemma 2}. {\it The action described above is free everywhere in
$Z_P$ outside of two submanifolds
$$
F_1= \left\{ \left( \begin{array}{rr}
  u_1 & v_1 \\
  u_2 & 0 \\
  u_3 & 0 \\
\end{array} \right) : |u_1|^2+|v_1|^2 = |u_2|^2=|u_3|^2=1 \right\},
$$
$$
F_2= \left\{ \left( \begin{array}{rr}
  0 & v_1 \\
  0 & v_2 \\
  u_3 & v_3 \\
\end{array} \right) : |v_1|^2 = |v_2|^2=|u_3|^2+|v_3|^2=1 \right\}.
$$
The points of the submanifolds $F_1$ and $F_2$ have the following
stabilizers:
$$
Fix(F_1)=\Gamma_1 = \{ (1, \pm (1,1)) \} = \Z_2, Fix(F_2) = \Gamma_2
= \{ (\pm (1,1),1)   \} = \Z_2.
$$ }

\vskip0.2cm

Thus, the factor-space $M'_3=Z_P/T^3$ is a quasitoric orbifold with
two singular submanifolds $F_1/T^3=S^2$ and $F_2/T^3=S^2$. Therefore,
the tubular neighborhoods of $F_1/T^3$ and $F_2/T^3$ are foliated with the fibres diffeomorphic to $\C^2/\Z_2$. If we resolve $M'_3$ along
the singular submanifolds (by blowing-up along the fibres of the tubular
neighborhoods), then we obtain a quasitoric manifold $N_3$ corresponding
to $N_3/T^3=Q_3$ and a principal $T^3$-bundle
$$
\pi_3: Z_{Q_3} \rightarrow N_3.
$$
As was already mentioned in the introduction, the manifold $Z_{Q_3}$
in this example is non-formal \cite{Baskakov,Buchshtaber-Panov}.

\vskip0.2cm

{\bf Theorem 1}. {\it The moment-angle manifolds $Z_{Q_1}$, $Z_{Q_2}$
and $Z_{Q_3}$ from Examples 1-3 admit $T^k$-invariant Riemannian
metrics of positive Ricci curvature ($k=11, 4, 5$, respectively)
such that the principal bundles $\pi_1$, $\pi_2$ and $\pi_3$ are
Riemannian submersions. }

\vskip0.2cm

The remaining part of the article is devoted to the proof of this theorem.

\section{Blowing-up at singular points of the manifolds with positive Ricci curvature}

Let us consider a  round  $(2n+1)$-sphere $S_R^{2n+1}$ of the radius
$R$ in the complex space $\C^n$. The standard Hermitian scalar
product induces a Riemannian metric $R^2 d s_{2n+1}^2$ on
$S_R^{2n+1}$. The circle $ S^1$ acts freely on the sphere
$S_R^{2n+1}$: for $ u \in S^1 we have: (z_0, \ldots, z_n) \mapsto (u
z_0, \ldots, u z_n)$. The factor-space under this action is the
complex projective space $\C P_R^n$ with the Fubini-Study metric
scaled by $R$. As usual, we will denote $\C P_1^n=\C P^n$.

Let us remind the (topological) construction of blowing-up the
complex projective space $\C P^n$. Consider the canonical complex
$\C$-bundle over $\C P^{n-1}$. It is well known that the
corresponding spherical sub-bundle is isomorphic to the Hopf
fibration $S^{2n-1} \rightarrow \C P^{n-1}$. Let $E$ be the total
space of the globular sub-bundle in the canonical bundle, that is,
$\partial E = S^{2n-1}$. We can cut off a geodesic ball $B$ of the
radius $\varepsilon$ in $\C P^n$ and glue $\C P^n \backslash B$ with
$E$ along the common boundary. The obtained manifold can be endowed
with the canonical smooth structure and is diffeomorphic to $\C P^n
\# \overline{\C P^n}$. If we now consider an orbifold $\C P^n/\Z_p$,
where the group $\Z_p$ acts with the isolated fixed point $p$, then
the previous construction can be repeated if we take a geodesic ball
with the center $p$. In this case $\partial ((\C P^n/\Z_p)
\backslash B) = S^{2n-1}/\Z_p$, so we should consider the $p$-th
tensor power of the canonical fibration $E^p$, $\partial(E^p) =
S^{2n-1}/\Z_p$ and glue $(\C P^n/\Z_p) \backslash B$ with $E^p$
along the common boundary. The resulted manifold $M$ models the
blowing-up of a singular point. Notice that there exists a
$T^n$-action on $\C P^n/\Z_p$ and $E^p$. Moreover, the restrictions
of this action to the identified boundaries are the same, so there
can be defined a natural $T^n$-action on the manifold $M$.

We will be interested in the following specific action of the group $\Z_n$ on $\C
P^n$. Let $\omega=e^{2\pi i/n}$, $\Z_n=\langle \omega \rangle$.
Consider the action on $\C^{n+1}$:
$$
\omega: (z_0, \ldots, z_n ) \mapsto (z_0, \omega z_1, \ldots, \omega
z_n).
$$
It is clear that this action induces an isometric action $\Z_n$ on $\C
P^n$, which is free outside of two subsets: $\bar{p}=[1:0:\ldots
: 0]$ and $\C P^{n-1}=\{ z_0=0\}$. Therefore, $\bar{p}$ is an
isolated fixed point, and we can perform the blow-up operation in its
neighborhood. The goal of the present section is to prove the following
theorem.

\vskip0.2cm

{\bf Theorem 2}. {\it There exists $\sigma>0$ such that for every
$\varepsilon>0$ the manifold $M$, obtained from $\C P^n/\Z_n$ by
blowing-up the singular point $\bar{p}$ as described above, admits a Riemannian metric $g$ with Ricci curvature bounded below by $\sigma$
and coinciding with the Fubini-Study metric outside of the geodesic
ball of radius $\varepsilon$ centered at the point $\bar{p}$.
Moreover, the torus $T^n$ acts on $M$ isometrically.

}

\vskip0.2cm

We begin the proof of Theorem 2 by the following lemma.

\vskip0.2cm

{\bf Lemma 3}. {\it The Fubini-Study metric on $\C P_R^n$ can be
expressed in the following way:
$$
g_R= dt^2 + R^2 \sin^2 \frac{t}{R} \ \cos^2 \frac{t}{R} \ d s_v^2 +
R^2 \sin^2 \frac{t}{R} \ d s_h^2,
$$
where $0 \leq t \leq R \pih$ is the distance to the point $\bar{p}$, whereas
$ds_v^2=d s_{2n-1}^2|_V$ and $d s_h^2=d s_{2n-1}^2|_H$ are the
restrictions of the standard sphere metric to the distribution of vertical
and horizontal tangent subspaces for the diagonal action of $S^1$ on the sphere
$S^{2n-1}$.

}

\vskip0.2cm

{\bf Proof of Lemma 3}. The tangent space to the sphere at the point $p
\in S_R^{2n+1}$ splits into the horizontal and vertical subspaces $T_p
S_R^{2n+1}=V_p \oplus H_p$ with respect to the Riemannian submersion
$\pi_n: S_R^{2n+1} \rightarrow \C P_R^n$:
$$
\begin{array}{c}
V_p=\{ t \cdot i p \ | \ t \in \R \}, \\
H_p= \{ u \in \C^{n+1} \ | \ \langle u, p \rangle_\C=0 \}.
\end{array}
$$

Let us fix the points $p_0=(R,0,\ldots, 0) \in S_R^{2n+1}$ and
$\bar{p}_0=\pi_n(p_0)$. For $0 \leq t \leq R \pi/2$ we put
$$
S_t= \{ (z_0, z_1, \ldots, z_n) | |z_0|=R \cos \frac{t}{R} \}
\subset S_R^{2n+1}.
$$
It is clear that $S_0$ is the $S^1$-orbit of the point $p_0$ and a
closed geodesic on the sphere; $S_{R \pi/2}$ is the equatorial sphere in
$S_R^{2n-1}$ under the imbedding $\C^{n} \subset \C^{n+1}$ as
a complex hyperplane $\{ z_0=0 \}$. For $0 < t < R \pi/2$   the
submanifold $S_t$ is the tubular hypersurface of radius $t$ around
$S_0$ (evidently, it coincides with the tubular hypersurface of
radius $R \pi/2 - t$ around $S_{R \pi/2}$), and it is isometric to the product
$S^1 \times S^{2n-1}$ with the metric
$$
g=R^2 ( \cos^2\frac{t}{R} \  ds_1^2 + \sin^2\frac{t}{R} \
ds_{2n-1}^2).
$$
Every shortest normal geodesic from $S_0$ to $S_t$ projects to a
geodesic in $\C P_R^n$ of the same length. Thus, $S_t$ projects to a
geodesic sphere $\bar{S}_t$ in $\C P_R^n$ of radius $t$ centered at
$\bar{p}_0$, which is the factor-space of $S_t$ under the action of
the group $S^1$. It is easy to see that, under the assumption $0<t<R
\pi/2$,   $\bar{S}_t$  is diffeomorphic to the sphere $S_R^{2n-1}$
with a  "skewed" metric. If we denote $d s_v^2=d s_{2n-1}^2|_V$ and
$d s_h^2=d s_{2n-1}^2|_H$ to be the restrictions of the standard
spherical metric to the distributions of vertical and horizontal
subspaces in $\bar{S}_t$, then it can be verified directly that the
metric on $\bar{S}_t$ looks like:
$$
\bar{g}= R^2 ( \sin^2 \frac{t}{R} \ \cos^2 \frac{t}{R} \  d s_v^2 +
\sin^2 \frac{t}{R} \ d s_h^2).
$$
So we obtain the following expression for the Fubini-Study metric on
$\C P_R^n$:
$$
g_R= dt^2 + R^2 \sin^2 \frac{t}{R} \ \cos^2 \frac{t}{R} \ d s_v^2 +
R^2 \sin^2 \frac{t}{R} \ d s_h^2,
$$
where $0 \leq t \leq R \pi/2$. The lemma is proved.

Further, we will consider metrics on $\R \times S^{2n-1}$ of the
following form:
$$
g=dt^2+ h(t)^2 ds_v^2+f(t)^2 ds_h^2\eqno{(2)}
$$
The following formulas are well known and can be checked easily.

\vskip0.2cm

{\bf Lemma 4}. {\it Let $X_0=\frac{\partial}{\partial t}$ be a
unit radial vector, let $X_1$ be a unit vertical vector (that is, $X_1
\in V$) and let $X_2$ be a unit horizontal vector ($X_2 \in H$). Then,
the Ricci curvature can be calculated by the following formulas.
$$
\begin{array}{c}
Ric(X_i, X_j)=0, \mbox{\ if \ } i \neq j, \\
Ric(X_0,X_0)=-\frac{h''}{h}-(2n-2)\frac{f''}{f}, \\
Ric(X_1,X_1)=-\frac{h''}{h}-(2n-2)\frac{f'h'}{fh}+(2n-2)\frac{h^2}{f^4},
\\
Ric(X_2,X_2)=-\frac{f''}{f} -\frac{f'h'}{fh}
-(2n-3)\frac{(f')^2}{f^2} +\frac{2n}{f^2}- 2\frac{h^2}{f^4}.\\
\end{array}\eqno{(3)}
$$

}

\vskip0.2cm

{\bf Proof of the Theorem 2}. First, we consider a metric $g$,
represented in the following special form:
$$
g=\frac{dr^2}{1-\phi(r)}+ r^2(1-\phi(r))ds_v^2+r^2ds_h^2,\eqno(4)
$$
where $\phi(r)$ is a smooth function. Then, the formulas (3) are transformed as follows:
$$
\begin{array}{c}
Ric(X_0,X_0)=Ric(X_1,X_1)=\frac{1}{2}(\phi''+(2n+1)\frac{\phi'}{r}),
\\
Ric(X_2,X_2)=\frac{\phi'}{r}+2n\frac{\phi}{r^2}.
\end{array}
$$
For instance, the Euclidean metric on $\C^n$ corresponds to $\phi(r)=0$,
whereas the metric $g_R$ on $\C P_R^n$, considered above, can be obtained by taking
$$
\phi(r)=\phi_R(r)=\frac{r^2}{R^2}, 0 \leq r \leq R.
$$
Note that $g_R$ is an Einstein metric with the cosmological constant
$\frac{2(n+1)}{R^2}$. Put $\psi(r)=r \phi'+2n \phi$. Then, the Ricci
tensor looks especially simple:
$$
\begin{array}{c}
Ric(X_0,X_0)=Ric(X_1,X_1)=\frac{\psi'}{2 r},
\\ Ric(X_2,X_2)=\frac{\psi}{r^2}.
\end{array}
$$
It is evident that the function $\psi_R(r)=2(n+1)\frac{r^2}{R^2}$
corresponds to the metric on $\C P_R^n$.

{\bf Remark 1}. If we put $\psi(r)\equiv 0$, then we obtain
a Ricci-flat metric
$$
g_n=\frac{dr^2}{1-\frac{1}{r^{2n}}}+r^2
\left(1-\frac{1}{r^{2n}}\right) ds_v^2 +r^2ds_h^2,
$$
which is nothing else but the $SU(n)$-holonomy metric discovered by Calabi  in \cite{Calabi}.
It was precisely the Calabi metric that brought us to the ansatz (4).

Put $r_1=\sqrt{R}$. Now, let us define a function $\psi_n(r)$ for $1 \leq r
\leq R$. Let
$$
\psi_n(r)=\psi_R(r), \mbox{\ if \ } r_1 \leq r.
$$
We have $\psi_n'(r_1)= 4(n+1)r_1/R^2>0$. Therefore, we can extend the
function $\psi_n$ smoothly into the segment $[1, r_1]$, so that it will
satisfy the following conditions in this segment:
$$
\begin{array}{l}
\psi_n'(r) >\frac{\kappa}{R^2},  \mbox{\ if \ }  r\geq 1, \\
\psi_n(r) \geq \kappa \frac{r^2-1}{R^2},  \mbox{\ if \ }  1 \leq r \leq r_1, \\
\psi_n(1)= 0, \\
\int_1^{r_1} \psi_n(s) s^{2n-1} ds = \eta,
\end{array}
$$
where $\kappa>0$ is independent on $R$, and $\eta$ is any prescribed
number satisfying the inequality
$$
\kappa \frac{r_1^{2n}}{R^2}
\left( \frac{r_1^2}{2n+2} - \frac{1}{2n} \right) < \eta <
\frac{\psi_n(r_1)}{2n} \left(r_1^{2n}-1\right).
$$
We would like to put
$$
\eta=\frac{r_1^{2(n+1)}}{R^2} -1.
$$
Indeed, the inequalities for $\eta$ are written in this case as:
$$
\kappa R^{n-1} \left( \frac{1}{2n+2} - \frac{1}{2n R}\right) <
R^{n-1}-1,
$$
$$
\left(R^{n-1} - 1\right) < \frac{n+1}{n} \left(R^{n-1} -
\frac{1}{R}\right).
$$
Thus, if $R>2$, then both inequalities are satisfied, and we can
construct a function $\psi_n$ with the required properties. Now we can put
$$
\phi_n(r)=\frac{1}{r^{2n}}\int_{1}^r \psi_n(s) s^{2n-1} ds +
\frac{1}{r^{2n}}.
$$
Obviously, the metric $g_n$, constructed on the basis of the last function, has nonnegative Ricci curvature, which is positive for $r>1$. Moreover, the Ricci curvature at the points $r=1$ vanishes only for the directions of the $X_2$ vectors. Now, let us verify that $\phi_n$ coincides with $\phi_R$ in the segment $r_1 \leq r \leq R$. Indeed, if $r \geq
r_1$, then
$$
\phi_n(r)=\frac{1}{r^{2n}}+ \frac{1}{r^{2n}}\int_{r_0}^{r_1}
\psi_n(s) s^{2n-1} ds + \frac{1}{r^{2n}}\int_{r_1}^r
2(n+1)\frac{s^{2n+1}}{R^2} ds =
$$
$$
=\frac{1}{r^{2n}}\left(1+\eta - \frac{r_1^{2(n+1)}}{R^2}\right)+
\frac{r^2}{R^2}=\frac{r^2}{R^2}.
$$

Thus, the metric $g_n$, corresponding to the function $\phi_n$, coincides with
the metric $g_R$ on $\C P_R^n$ for $r \geq r_1$. Note that
$$
\phi_n (1)=1, \phi_n^\prime (1) = -2 n.
$$

To get rid of the effect of the Ricci curvature vanishing for $r=1$, we consider
a modified metric:
$$
g_n'=\frac{dr^2}{1-\phi_n (r)}+ r^2(1-\phi_n
(r))ds_v^2+(1-\delta(r)) r^2ds_h^2,
$$
where $\delta (r)\geq 0$ is a smooth function. Choose a
sufficiently small $\nu >0$. Take a non-increasing function
$\delta_\nu (r)$, $1 \leq r \leq R$, such that the following
conditions are satisfied:
$$
\begin{array}{l}
\delta_\nu (r) = 0 \mbox{ \ if \ } r \geq r_1, \\
\delta_\nu (r) >0 \mbox{ \ if \ } 1 \leq r < r_1, \\
\delta_\nu (r) \mbox{\ is \ constant } \mbox{ for \ } 1 \leq  r \leq 2, \\
\|\delta_\nu (r)\|_{C^\infty} \leq \nu.
\end{array}
$$
It is clear that if we choose $R>4$, then for each $\nu>0$ there exists a function
$\delta_\nu (r)$ with the above properties. Then (3) shows that
for $1\leq r \leq 2$
$$
Ric' (X_2, X_2) = \frac{\psi_n }{r^2} + \frac{2 \delta_\nu (1)
}{(1-\delta_\nu (1)) r^2} \left( n - \frac{2-\delta_\nu
(1)}{1-\delta_\nu (1)} \left(1-\phi_n \right)\right).
$$
Since the function $\phi_n$ is strictly positive everywhere on $[1,R]$,
then we can assume $\phi_n \geq c(R)>0$. Consequently, if we take
$\nu>0$ such that
$$
\delta_\nu (r) \leq \delta_\nu (1) \leq \frac{2 c(R)}{1+c(R)},
$$
then it can be immediately checked that for $1 \leq r \leq 2$
$$
Ric' (X_2, X_2) \geq \frac{\psi_n }{r^2} + \frac{2 \delta_\nu (1)
}{(1-\delta_\nu (1)) r^2} \left( n - 2 \right) \geq \frac{2
\delta_\nu (1) }{(1-\delta_\nu (1))} \frac{ n - 2}{4} >0.\eqno{(5)}
$$
Then, for  $2 \leq r \leq R$ the curvature of $g_n$ is
$$
Ric(X_2, X_2)= \frac{\psi_n }{r^2} \geq \kappa \frac{r^2-1}{R^2}
\geq \frac{3\kappa}{R^2}
$$
By continuity, for all sufficiently small $\nu>0$ (let us remind that
$\nu$ is independent of $R$) the following estimate from below holds for the Ricci
curvature of the metric $g_n'$  for $2\leq r \leq R$:
$$
Ric'(X_2, X_2) \geq \frac{2\kappa}{R^2}\eqno{(6)}
$$
On the other hand, since the Ricci curvatures $Ric(X_0,X_0) = Ric(X_1,
X_1)$ of the metric $g_n$ are bounded below by $\frac{\kappa}{R^2}$, then
for all sufficiently small $\nu$ we have analogous estimates for the
metric $g_n'$:
$$
Ric'(X_0,X_0) = Ric'(X_1, X_1) \geq \frac{\kappa}{2 R^2}.\eqno{(7)}
$$
Thus, if we fix a sufficiently small number $\nu$, then by inequalities
(5), (6), (7) we obtain a common estimate from below for the Ricci curvature:
$$
Ric'(X,X) \geq \frac{\sigma}{R^2},\eqno{(8)}
$$
where $|X|=1$, and the constant $\sigma>0$ is chosen to be independent of $R>4$. In
particular,
$$
\frac{2 \delta_\nu (1) }{(1-\delta_\nu (1))} \frac{ n - 2}{4} \geq
\sigma.\eqno{(9)}
$$
Now, let us discuss the smoothness of the metric $g_n'$ in the neighborhood of
$r=1$. Firstly, since $h$ vanishes at the point $r=1$ and $f$ is
strictly positive, the integral circles of the vertical distribution $V$
(that is, the fibres of the submersion) collapse to a point when $r
\rightarrow 1$; and the metric $g_n'$ is well defined on the blow-up of the
space $\C P_R^n$ at the point $\bar{p}$. The smoothness of such
metrics was investigated, for example, in \cite{Bazaikin}. The
criterion of the smoothness is as follows:
$$
\frac{d}{dt}\big|_{r=1} f(t)=0, \frac{d}{dt}\big|_{r=1} h(t)= \pm 1
$$
(strictly speaking, this conditions imply only the $C^1$-smoothess, but the
quasi-linearity of the Ricci curvature operator enables us to $\C^{\infty}$-smooth such metric while preserving an estimate similar to (8)).
Direct calculations show that
$$
\frac{d}{dt}\big|_{r=1} f(t)=0,
$$
$$
\frac{d}{dt}\big|_{r=1} h(t)= - \frac{\phi'(1)}{2} = n.
$$
We can see that to ensure smoothness we have to contract the vertical circles of the submersion by $n$ times; that is, the metric $g_n'$ is a smooth metric with
positive Ricci curvature on the space obtained from $\C P_R^n/\Z_n$ by resolving the singularity.

Now, rescale the metric $g_n'$ by $1/R^2$. We obtain a Riemannian metric
$g$ on $\C P^n/\Z_n$ with the blown-up singular point and, moreover, the blowing-up is performed inside a neighborhood of radius
$$
\frac{r_1}{R}=\frac{1}{\sqrt{R}}=\varepsilon \rightarrow 0,
$$
where $R \rightarrow \infty$. The estimate (8) turns into the following estimate
$$
Ric(X,X) \geq \sigma,
$$
where $\sigma$ does not depend on $\varepsilon$. It remains to
note that we have been preserving the $T^n$-invariance of the metric during all our
operations with the metric. The theorem is proved.

{\bf Remark 2}. Using (9), we can estimate
$$
f(1) =\sqrt{1-\delta_\nu (1)} \leq \sqrt{\frac{n-2}{n-2+2 \sigma}}.
$$
After the reduction of the metric by $R^2$ times, we have the
following estimate of the "size" $D$ of the space $\C P^{n-1}$,
which is pasted in instead of the singular point:
$$
D \leq \frac{1}{R}\sqrt{\frac{n-2}{n-2+2 \sigma}}  = \varepsilon^2
\sqrt{\frac{n-2}{n-2+2 \sigma}}\eqno{(10)}
$$
This estimate will be useful for us in the sequel.

\section[]{The construction of the Ricci-positive Riemannian metrics}

We will need several previously known results. Firstly, we need to lift
a Riemannian metric of positive Ricci curvature from the base to the total
space of a Riemannian submersion while preserving the Ricci-positivity.
Similar constructions were considered in the papers \cite{Cheeger, Poor,
Berard-Bergery, Nash, Gilkey-Park-Tuschmann}. In
\cite{Berard-Bergery} the possibility of such lift was proved for
the principal torus bundles (the only bundles we will need in our paper), but
without any statement concerning the invariance of the resulting metric. Therefore, to
obtain the invariant metrics on the moment-angle manifolds, we will use the
following stronger result.

\vskip0.2cm

{\bf Theorem 3} \cite{Gilkey-Park-Tuschmann}. {\it Let $(Y,g_Y)$ be
a compact connected Riemannian manifold with positive Ricci
curvature. Let $P$ be a principal bundle over $Y$ with the compact
connected structure group $G$ so that $\pi_1(P)$ is finite. Then,
there exists a $G$-invariant metric $g_P$ on $P$ such that $g_P$ has
positive Ricci curvature and $\pi: (P,g_P) \rightarrow (Y,
g_Y)$ is a Riemannian submersion. }

\vskip0.2cm

The following theorem allows changing arbitrary Riemannian metric in
a small neighborhood of the manifold to the some model metric.

\vskip0.2cm

{\bf Theorem 4} \cite{Gao}. {\it In the ball $U(0,\rho_0) = \{ x \in
\R^n| |x| \leq \rho_0 \}$, there are considered two Riemannian metrics $g_0$ and
$g_1$ with positive Ricci curvature and with the same $1$-jets
$J^1(g_0)$ and $J^1(g_1)$ at the point $0$. Then there exists a
Riemannian metric $\bar{g}$ in $U(0,\rho_0)$ with positive Ricci
curvature and  $0<\rho_2 < \rho_1 < \rho_0$, such that $\bar{g}=g_1$
for $|x| < \rho_2$ and $\bar{g}=g_0$ for $|x| > \rho_1$.

}

\vskip0.2cm

{\bf Remark 3}. In \cite{Gao} Theorem 4 was proved for negative
Ricci curvature, but the proof is valid for positive Ricci
curvature. In addition, we will need to know slightly more about the construction used in the proof of Theorem 4. Firstly, the metric
$\bar{g}$ has the form:
$$
\bar{g}=(1-s)g_0 + s g_1,
$$
where $s=\psi(|x|)$ for some smooth function $\psi : \R \rightarrow
[0,1]$. This immediately implies that if the metrics $g_0, g_1$ are
invariant with respect to the action of a group $G$ preserving $|x|$, then
$\bar{g}$ will be also invariant with respect to $G$. Secondly,
the coefficients of the Levi-Chivita connection and the Ricci tensor coefficients of the metrics $g_1, g_2, \bar{g}$ can be chosen to be arbitrary close to each other,
independently of $\rho_0$.

Now, let us examine one by one the examples from Section 2.

{\bf Examples 1 and 2}. As it was described in Section 2, the orbifolds $M'_1$
and $M'_2$ have isolated singular points with the neighborhoods of the
form $\C^3/\Z_3$. Applying Theorem 4 (to the $\Z_3$-covering of the small
neighborhood of every singular point), we can assume that the neighborhood of each
singular point is isometric to a geodesic ball in $\C
P^3/\Z_3$. Applying Theorem 2 yields some Ricci-positive invariant
metric on $N_1$ and $N_2$, which are obtained by blowing-up the singular points. Finally, Theorem 3 provides us with the metrics needed on $Z_{Q_1}$ and $Z_{Q_2}$.

{\bf Example 3}. Take $\eps > 0$. Construct a smooth function $f(t)$
on the segment $[0, \pih]$, satisfying the following conditions:
$$
\begin{array}{l}
f(t) = 1 \mbox{\ for \ }  0 \leq t \leq \varepsilon, \\
f(t)=\cos t \mbox{\ for \ } \pih-\eps \leq t \leq \pih, \\
f'(t) < 0, f''(t) < 0 \mbox{\ for \ } \varepsilon < t <
\pih-\varepsilon.
\end{array}
$$
It is clear that such a function exists for all sufficiently small
$\eps$. Consider the following Riemannian metric on the space $[0,
\pih]^3 \times (S^1)^6$ with the coordinates $(t_1,$ $t_2$, $t_3$,
$\phi_1, \psi_1, \phi_2, \psi_2, \phi_3, \psi_3)$:
$$
\begin{array}{rl}
g = & dt_1^2 + \cos (t)^2 d\phi_1^2 + f(\pih - t_1)^2 d\psi_1^2 +  \\
    & + dt_2^2 + f(t_2)^2 d\phi_2^2 + f(\pih - t_2)^2 d\psi_2^2 +  \\
    & + dt_3^2 + f(t_3)^2 d\phi_3^2 + \sin(t_3)^2 d\psi_3^2. \\
\end{array}
$$
It is easy to see that for each $i=1, 2, 3$ the coordinates $(t_i,
\phi_i, \psi_i)$ are well defined on the corresponding spheres $S^3$; $g$
is a smooth metric on $Z_P = S^3 \times S^3 \times S^3$ and has
nonnegative sectional curvature because of the restrictions imposed on
$f(t)$.

In our coordinates, the subsets $F_1$ and $F_2$ are specified by the
equations $t_2 = t_3 = 0$ and $t_1 = t_2 = \pih$, respectively.

{\bf Remark 4}. It is clear that the transformation $(t_1, t_2, t_3,
\phi_1, \psi_1, \phi_2, \psi_2, \phi_3, \psi_3) \mapsto (\pih-t_3,
\pih - t_2, \pih - t_1, \psi_3, \phi_3, \psi_2, \phi_2, \psi_1,
\phi_1)$ is an isometry of $Z_P$, which interchanges $F_1$ and
$F_2$. Therefore, it suffices to perform further reasoning (due to
its locality) only for $F_1$.

Define a neighborhood $U \subset Z_P$ of the set $F_1$ by the relations
$t_2 < \eps$ and $t_3 < \eps$. Obviously, $U$ is isometric to
the direct product $S^3 \times S^1 \times S^1 \times D^2 \times D^2$
with the metric
$$
\begin{array}{rl}
g|_{U} = & \left( dt_1^2 + \cos (t)^2 d\phi_1^2 + f(\pih - t_1)^2
                     d\psi_1^2 \right) + \left(  d\phi_2^2 +  d\phi_3^2 \right) + \\
    & + \left( dt_2^2 + \sin (t_2)^2 d\psi_2^2 \right) + \left( dt_3^2 + \sin(t_3)^2 d\psi_3^2 \right) \\
\end{array}
$$
(that is, every disk $D^2$ is isometric to a geodesic disk on the
two-dimensional unit sphere). Put $S=D^2 \times D^2$ with the induced metric. In what follows, we will also consider in $S$ the coordinates $(t,
\theta, \psi, \phi)$ given by $t_2= t \cos(\frac{\theta}{2}), t_3 =
\sin(\frac{\theta}{2}), \psi_2 = \frac{\psi+\phi}{2}, \psi_3 =
\frac{\psi-\phi}{2}$.

It is clear that $S$ has nonnegative sectional curvature and
strictly positive Ricci curvature (to be precise, the Ricci
curvature of $S$ equals one). Applying Theorem 4, we can deform the
metric on $S$ in such a way that, preserving the Ricci-positivity,
it will not change outside the $\varepsilon/2$-neighborhood of the
point $p = \{ t=0 \} \in S$ and will be isometric to a geodesic ball
of radius $\varepsilon/4$ in $\C P^2$ inside the
$\varepsilon/4$-neighborhood. According to Theorem 2, we can now
blow-up $S/\Z_2$ inside the $\varepsilon/4$-neighborhood preserving
Ricci-positivity, obtaining a manifold $S'$. Now we replace each
"direct" factor $S$ in $U$ by a double-covering over $S'$ (branched
along the vertical circles of the Hopf fibration) and obtain $U' =
S^3 \times S^1 \times S^1 \times S'$. Since the Riemannian manifolds
$S$ and the double-covering over $S'$ are isometric in the
$\varepsilon/2$-neighborhood, we can paste in $U'$ instead of $U$ in
$Z_P$, obtaining $Z_P'$. Strictly speaking, $U'$ is a
double-covering, branched over $F_1$, but after factorizing with
respect to the action (1) of the group $T^3$ we get a smooth
manifold. Let us remark that this whole construction can be
performed for any arbitrary small $\varepsilon$.

\vskip0.2cm

{\bf Lemma 4}. {\it The factor-space $\tilde{U}=U'/T^3$ with the
induced metric has strictly positive Ricci curvature for a sufficiently small $\varepsilon$. }

\vskip0.2cm

{\bf Proof}. The metric induced on the space $\tilde{U}$ is
characterized by such property that the factorization map $U' \rightarrow
\tilde{U}$ is a Riemannian submersion. We use the following analogue
of O'Neill’s formula for the Ricci tensor \cite{Besse}:
$$
\tilde{R}(\tilde{X},\tilde{X})=R(X,X)+2 (AX, AX)+(TX,TX)-(D_X N, X).
$$
Here $\tilde{R}$ and $R$ are the Ricci tensors of $\tilde{U}$ and $U'$;
$A$ and $T$ are the fundamental tensors connected with the Riemannian submersion;
$\tilde{X}$ is a tangent vector field on $\tilde{U}$ with the horizontal
lifting $X$ on $U'$; $N$ is the mean curvature vector of the submersion
fibers, that is, $N=\sum_{i} T_{U_i} U_i$, where $U_i$ is an orthonormal
basis in the vertical subspace of submersion. As a (non-orthonormal) basis of the vertical vector fields, we can choose the following one:
$$
V_1=\partial_{\phi_1}+ 2
\partial_{\psi_1} + \partial_{\phi_2},
$$
$$
V_2= \partial_{\psi_1} +\partial_{\phi_2} +  \partial_{\phi_3} +
\partial_{\psi_2},
$$
$$
V_3= -\partial_{\psi_1} - \partial_{\phi_2} +
\partial_{\phi_3} + \partial_{\psi_3}.
$$
Thus,
$$
\tilde{R}(\tilde{X},\tilde{X}) \geq R(X,X) - (D_X N, X).\eqno{(11)}
$$
Note first that $Z_P$ has nonnegative Ricci curvature and
$R(X,X)=0$ exactly in two cases: 1) $t_2 \leq \varepsilon, t_3 \leq
\varepsilon$ and $X=\alpha
\partial_{\phi_2} + \beta
\partial_{\phi_3}$; 2) $ t_1 \geq \pih-\varepsilon, t_2 \geq \pih -
\varepsilon$ and $X=\alpha \partial_{\psi_1} + \beta
\partial_{\psi_2}$, where $\alpha, \beta$ are some coefficients.
One can verify immediately that the non-vanishing vectors $X$ from 1), 2)
cannot be horizontal. Consequently, there exists a constant $\sigma>0$ such that the Ricci curvature of the space $Z_P$ in the horizontal direction
is no less than $\sigma$. Therefore, in view of Remark 3 and after
applying Theorem 4 to deform the metric, we can achieve that the Ricci
curvature of the resulting metric in horizontal direction is no less
than $\sigma/2$ independently of how small $\varepsilon$ is. At last, applying
Theorem 2, the Ricci curvature in horizontal direction after blowing-up is
again estimated from below by some constant $\kappa>0$ independent of
$\varepsilon$. Thus, we can estimate from below the first summand in the
right-hand side of (11):
$$
R(X,X) \geq \kappa.\eqno{(12)}
$$

Now, suppose that we choose some orthonormal basis in the space of
vertical vectors:
$$
V'_i= Y_i + Z_i,
$$
where $Y_i$, $Z_i$ are the vector fields tangent to $S^3 \times T^2$
and $S'$, respectively, where $i=1, 2, 3$. Then
$$
N = \sum_{i=1}^3 \left( T_{Y_i} Y_i + T_{Z_i} Z_i \right) = N_1 +
N_2
$$
is the corresponding decomposition of the mean curvature vector. We see
that the quadratic form $R(X,X) - (D_X N, X)$ splits into three blocks,
corresponding to the tangent subspaces of $S^3$, $T^2$ (flat torus) and
$S$.

If we apply the construction of Theorem 4, then, in view of Remark
3, we obtain a metric $\bar{g}$ on $S$ which is a "linear
interpolation"  of the metrics $g_0$ and $g_1$. Therefore, in our
case it will have the form
$$
\bar{g} = \left( dt^2 + f(t)^2 d \theta^2 \right) + \sum_{i,j=1}^2
g_{ij} (t,\theta) d \xi^i d \xi^j,
$$
where $\psi=\xi^1, \phi=\xi^2$. A simple calculation shows that
$$
D_{\xi^i} \xi^j = -\frac{1}{2} \mbox{ grad} (g_{ij})
$$
(the gradient is taken with respect to the metric $dt^2+h^2 d\theta^2$). Let us put
$$
Z_i= \sum_{j=1}^2 \alpha_i^j \partial_{\xi^j}.
$$
Since $T^3$ acts with a fixed point at the pole $p \in  S$, the vector fields $Z_i$
vanish at $t=0$. Thus $\alpha_i^j (p) = 0$ for all $i, j$. Then,
$$
N_2 = - \frac{1}{2} \sum_{i=1}^3 \sum_{j,k=1}^2 \alpha_i^j
\alpha_i^k \mbox{ grad} (g_{jk}).
$$
Therefore,
$$
- (D_X N_2, X)|_{t=0} = \left( - X (N_2, X) + (N_2, D_X
X)\right)_{t=0} =
$$
$$
=  \frac{1}{2} \sum_{i=1}^3 \sum_{j,k=1}^2 \left( 2 \alpha_i^j
g_{jk} X (\alpha_i^k)  + \alpha_i^j \alpha_i^k X
(g_{jk})\right)_{t=0} = 0.
$$
This means that, taking a sufficiently small $\varepsilon>0$, we can
achieve that throughout the domain $U$ the following inequality for the metric $\bar{g}$

holds:
$$
|(D_X N_2, X)| \leq \frac{\kappa}{2}.\eqno{(13)}
$$

Now let us see what happens with the quantity $(D_X N_2, X)$ during the blow-up process as described in Theorem 2.
The deformation does not take the metric out of the following class of metrics:
$$
\bar{g}=dt^2 + \frac{1}{4} h(t)^2 \left( d \psi + \cos (\theta)
d\phi \right)^2 + \frac{1}{4} f(t)^2 \left( d \theta^2 +
\sin(\theta)^2 d\phi^2 \right)
$$
(this metric is exactly the metric (2) under the assumption $n=2$,
expressed in the coordinates $(t,\theta, \psi, \phi)$). Recall that here
$h(0)=0$, $f'(0)=0$, and $f(0)=f_0>0$ can be taken arbitrary small
while decreasing $\varepsilon$ (see (10) and Remark 2). On the other
hand, the coordinates of the field $V'_i$ corresponding to $\partial_{\psi}$ vanish
at $t=0$, because the sphere $\{ t=0 \}$ is fixed by the action of
$\partial_{\psi}$. Thus, $\alpha_i^1 (p) =0$. Now, as the radius of the sphere
$t=0$ equals $f_0$, we can make this sphere arbitrary small as well; and so $\alpha_i^2$ at $t=0$ can be chosen arbitrary small while decreasing $\varepsilon$. So, in this case we
have
$$
- (D_X N_2, X)|_{t=0} = \left( - X (N_2, X) + (N_2, D_X
X)\right)_{t=0} =
$$
$$
=  \frac{1}{2} \sum_{i=1}^3 \sum_{j,k=1}^2 \left( 2 \alpha_i^j
g_{jk} X (\alpha_i^k)  + \alpha_i^j \alpha_i^k X (g_{ij}) +
\alpha_i^j \alpha_i^k \left( \mbox{grad} (g_{jk}), D_X X
\right)\right)_{t=0}
$$
Since it is guaranteed that all the derivatives of the functions $f(t), h(t)$ are bounded independently of $\varepsilon$, all the derivatives of the coefficients of
$g_{ij}$ are bounded along the unit vector fields. This means that by diminishing $\varepsilon$ we can get the estimate (13) in this case as well.

Now, consider the component $N_1$. Take the following basis of the vertical
vectors:
$$
V_1=2 \partial_{\phi_3}+\partial_{\psi_2} +
\partial_{\psi_3}.
$$
$$
V_2= \partial_{\psi_1}+\partial_{\phi_2}+\partial_{\phi_3} +
\partial_{\psi_2}.
$$
$$
V_3=\partial_{\phi_1}-\partial_{\phi_2}-2\partial_{\phi_3} -2
\partial_{\psi_2},
$$
By the standard orthogonalization method we pass to an orthonormal basis at each point of
the vertical vector fields $V_1', V_2', V_3'$:
$$
V_1'=Z'_1,
$$
$$
V_2'=y_1 \partial_{\psi_1}+ Z'_2,
$$
$$
V_3'=z_1 \partial_{\phi_1} +z_2 \partial_{\psi_1}+ Z'_3,
$$
where the fields $Z'_i$ are tangent to $T^2 \times S'$, and the coefficients
$y_1, z_1, z_2$ are the functions of $(t_1, t_2, t_3)$. Consequently,
$$
N_1=  \left(z_1^2 \sin(t_1) \cos (t_1)+w^2 f\left(\pih-t_1\right)
f'\left(\pih-t_1\right)\right)
\partial_{t_1}
$$
for some $w$, which can be explicitly expressed via $y_1, z_1, z_2$.
As above, we see that $(D_X N_1, X) \neq 0$ if only $X=\partial_{t_1}$. So,
we need only to check the Ricci-positivity for $X=\partial_{t_1}$. However, for all the horizontal vectors $Y$, which are orthogonal to
$\partial_{t_1}$, the sectional curvature $K(\partial_{t_1}, Y) \geq 0$,
and, moreover, $K(\partial_{t_1},
\partial_{\phi_1})=1$. Thus, by O'Neill’s formula we have
$$
\tilde{R}(\tilde{\partial}_{t_1}, \tilde{\partial}_{t_1}) =
\sum_{\tilde{Y}\perp \tilde{\partial}_{t_1}}
K(\tilde{\partial}_{t_1}, \tilde{Y})  \geq \sum_{Y\perp
\partial_{t_1}, V_1, V_2, V_3} K(\partial_{t_1}, Y)  > 0,\eqno{(14)}
$$
because there exist horizontal vectors $Y$ with an everywhere non-vanishing
coordinate with respect to $\partial_{\phi_1}$.

To complete the proof of the lemma, we consider an arbitrary horizontal vector
$$
X= \alpha_1 \partial_{t_1} +\alpha_2
\partial_{\phi_1}+ \alpha_3 \partial_{\psi_1} +
\alpha_4 \partial_{\phi_2} +\alpha_5 \partial_{\phi_3} + X_2,
$$
where $X_2$ is tangent to $S'$. Then, using (11), (12), (13) and (14),
we get
$$
\tilde{R}(\tilde{X},\tilde{X}) \geq \alpha_1^2
\tilde{R}(\tilde{\partial}_{t_1},\tilde{\partial}_{t_1})+
$$
$$
+
 R(\alpha_2
\partial_{\phi_1}+ \alpha_3 \partial_{\psi_1} +
\alpha_4 \partial_{\phi_2}, \alpha_2
\partial_{\phi_1}+ \alpha_3 \partial_{\psi_1} +
\alpha_4 \partial_{\phi_2}) + |\tilde{X_2}|^2 \frac{\kappa}{2}.
$$
It is clear that if $\alpha_1=0$ and $X_2=0$, then  $\alpha_2
\partial_{\phi_1}+ \alpha_3 \partial_{\psi_1} +
\alpha_4 \partial_{\phi_2}$ is a horizontal vector; and, therefore,
the middle summand in the right-hand side of the last inequality is
bounded below by $\sigma$. If $\alpha_1$ (respectively $|X_2|$) does
not vanish, then the middle summand is at least nonnegative and the
first (respectively, the third) summand is strictly positive. The
lemma is proved.

If we perform the blow-up construction described above for the neighborhood of $F_2$ as well, we obtain a quasitoric manifold $N_3$ with positive Ricci curvature.
Application of Theorem 3 finishes the proof of Theorem 1.

{\bf Remark 5}. The authors have been able to construct many other
actions of $T^3$ on $Z_P$, corresponding to cutting-off different
sets of vertices of the cube. It is clear that in all these cases, by utilizing similar reasoning, one can succeed in constructing the Ricci-positive Riemannian metrics on the corresponding moment-angle manifolds.

\end{document}